# Probability Bracket Notation: the Unified Expressions of Conditional Expectation and Conditional Probability in Quantum Modeling

Xing M. Wang
Sherman Visual Lab, Sunnyvale, CA, USA

## Table of Contents

## Abstract

After a brief introduction to *Probability Bracket Notation* (*PBN*), *indicator operator* and *conditional density operator* (*CDO*), we investigate probability spaces associated with various quantum systems: system with one observable (discrete or continuous), system with two commutative observables (independent or dependent) and a system of indistinguishable non-interacting many-particles. In each case, we derive unified expressions of *conditional expectation* (*CE*), *conditional probability* (*CP*), and *absolute probability* (*AP*): they have the same format for discrete or continuous spectrum; they are defined in both Hilbert space (using Dirac notation) and probability space (using *PBN*); and they may be useful to deal with *CE* of non-commutative observables.

## 1. Introduction to *PBN*, indicator operator and *CDO*

In our previous papers [1-5], we proposed and studied the applications of *Probability Bracket Notation* (*PBN*), a new set of symbols for probability modeling, inspired by *Dirac notation* in Quantum Mechanics (QM) [6-7]. In this paper, we will use both notations to derive unified expressions for absolute probability (*AP*), conditional probability (*CP*), and *conditional expectation* (*CE*) of observables for various quantum systems.

In this section, to prepare our discussion, we will give a brief review of the core *PBN* formulas for *stable* (time-independent) probability space, and then we will introduce the *indicator operator* and the *conditional density operator* (*CDO*), defined in both Hilbert space (using Dirac notation) and in probability space (using *PBN*).

***Introduction to PBN***: Let random variable ($R.V.$) $X$ be an observable associated with a stable probability space $(\Omega, \hat{X}, P)$, we have proposed [1-3]:

1. ***Conditional Probability*** of event $A$ under evidence $B$ is identified as a ***probability bracket*** from a ***P-bra*** and a ***P-ket*** :

$$P-bra: P(A| \qquad P-ket: |B) \qquad P-bracket: P(A|B) \equiv (A|B) \tag{1.1a}$$
$$P(A|B) = 0 \;\; if \;\; A \cap B = \varnothing$$
$$Discrete\ R.V. \quad P(A|B) = 1 \;\; \text{if} \;\; A \supset B \neq \varnothing$$
$$Continuous\ R.V. \quad P(A|B) = 1 \;\; \text{if} \;\; A \supset \bar{a} \supset B \neq \varnothing \;\; \& \int_{\bar{a}} dx > 0$$
$$P(A|B) = P(A|\hat{I}|B), \quad \text{where } \hat{I} \text{ is an unit operator} \tag{1.1b}$$

2. ***The basis*** (*orthonormality*):

$$Discrete\ R.V. \quad X|x_i) = x_i|x_i); \quad P(x_i|x_j) = \delta_{i,j} \tag{1.2a}$$
$$Continuous\ R.V. \quad X|x) = x|x); \quad P(x|x') = \delta(x-x') \tag{1.2b}$$

3. ***The unit operator*** (*completeness*) :

$$Discrete\ R.V. \qquad I = \sum_{x_i \in \Omega} |x_i)P(x_i| \tag{1.3a}$$
$$Continuous\ R.V. \quad I = \int_{x \in \Omega} dx\, |x)P(x| \tag{1.3b}$$

4. ***Expectation value*** of $R.V.\ X$ (in discrete and continuous cases):

$$E[X] \equiv \langle X \rangle \equiv P(\Omega|X|\Omega) = \sum_{x_i \in \Omega} P(\Omega|X|x_i)P(x_i|\Omega) = \sum_{x_i \in \Omega} x_i\, P(x_i) \tag{1.4a}$$
$$E[X] \equiv P(\Omega|X|\Omega) = \int_{x \in \Omega} dx\, P(\Omega|X|x)P(x|\Omega) = \int_{x \in \Omega} dx\, x\, P(x) \tag{1.4b}$$

5. ***Conditional expectation*** (*CE*) value of $X$ under evidence $A$:

$$E[X|A] \equiv P(\Omega|X|A) \tag{1.5}$$

.

6. ***Equivalence of probability distribution function***: The probability distribution function (*PDF*) [6-7]can be expressed in both probability space and in related Hilbert space as follows [3]:

$$P(x_i|\Omega) = P(x_i) = |\langle x_i|\Psi\rangle|^2 = |\Psi(x_i)|^2, \qquad \text{discrete spectrum} \tag{1.6a}$$
$$P(x|\Omega) = P(x) = |\langle x|\Psi\rangle|^2 = |\Psi(x)|^2, \qquad \text{cotinuous spectrum} \tag{1.6b}$$

To derive our unified *CE* expression, we need to introduce more definitions and expressions in *PBN*. In Ref [4], we have discussed *CE* in probability space using *PBN*, and introduced the *indicator operator,* based on the *indicator function* used in literature of probability theories [8]:

$$\boldsymbol{1}_A = \begin{cases} \sum_{x_i \in A} | x_i ) P(x_i | = \sum_{x_i \in A} \boldsymbol{1}_{x_i}, & \text{discrete spectrum} \quad (1.7a) \\ \int_{x \in A} dx \, | x ) P(x | = \int_{x \in A} dx \, \boldsymbol{1}_x, & \text{continuous spectrum} \quad (1.7b) \end{cases}$$

With the help of indicator operator, the *CE* of observable $\hat{X}$ given $A$ can be expressed as (see Eq. (3.1.3b) of [4] or Eq. (3.2.6) in [8]):

$$P(\Omega | X | A) = \frac{P(\Omega | X \boldsymbol{1}_A | \Omega)}{P(A | \Omega)}, \quad \text{where } P(A | \Omega) > 0 \quad (1.8)$$

Mapping (1.8) to induced Hilbert space [3], we obtain the indicator operator defined in Hilbert space using Dirac notation:

$$\boldsymbol{1}_A = \begin{cases} \sum_{i \in A} | x_i \rangle \langle x_i | = \sum_{i \in A} \boldsymbol{1}_i, & \text{discrete spectrum} \quad (1.9a) \\ \int_{x \in A} dx \, | x \rangle \langle x | = \int_{x \in A} dx \, \boldsymbol{1}_x, & \text{continuous spectrum} \quad (1.9b) \end{cases}$$

***Indicator operator*** is a natural extension of the *unit operator*, which now can be written as a special case of indicator operator:

$$\hat{I} = \boldsymbol{1}_\Omega = \begin{cases} \sum_{x_i \in \Omega} | x_i \rangle \langle x_i | = \sum_{x_i \in \Omega} \boldsymbol{1}_i, & \text{discrete spectrum} \quad (1.10a) \\ \int_{x \in \Omega} dx \, | x \rangle \langle x | = \int_{x \in \Omega} dx \, \boldsymbol{1}_x, & \text{continuous spectrum} \quad (1.10b) \end{cases}$$

The ***orthogonal projectors*** in (1.9-10) have the following general properties, valid in both Hilbert and probability spaces:

$$\boldsymbol{1}_i \boldsymbol{1}_k = \delta_{ik} \boldsymbol{1}_k, \quad \boldsymbol{1}_i \boldsymbol{1}_i = \boldsymbol{1}_i, \quad \text{discrete spectrum} \quad (1.11a)$$

$$\boldsymbol{1}_x \boldsymbol{1}_{x'} = \delta(x - x') \boldsymbol{1}_{x'}, \quad \boldsymbol{1}_x \boldsymbol{1}_x = \delta(0) \boldsymbol{1}_i, \quad \text{continuous spectrum} \quad (1.11b)$$

As seen in (1.11a), the discrete orthogonal projector $\boldsymbol{1}_i$ is self-adjoint. But, as seen in (1.11b), the continuous orthogonal projector $\boldsymbol{1}_x$ is not.

If set $A$ and $B$ are sets of orthogonal projectors as in (1.9), then we have:

Discrete: $\boldsymbol{1}_A \boldsymbol{1}_B = \sum_{x_i \in A} | x_i \rangle \langle x_i | \sum_{x_k \in B} | x_k \rangle \langle x_k | = \sum_{\substack{x_i \in A \\ x_k \in B}} | x_i \rangle \delta_{ik} \langle x_k | = \boldsymbol{1}_{A \cap B}$ (1.12a)

Continuous: $\boldsymbol{1}_A \boldsymbol{1}_B = \int_{x \in A \cap B} dx \, | x \rangle \langle x | = \boldsymbol{1}_{A \cap B}$ (1.12b)

***Density operator*** is used for the expression of expectation value of observables in Hilbert space (§22 of [6] and §11 of [7]). Using Dirac notation, it reads:

In Hilbert space: $\rho = | \Psi \rangle \langle \Psi |, \quad \mathrm{Tr}[\rho^2] = \mathrm{Tr}\rho = \langle \Psi | \Psi \rangle = 1$ (1.13a)

Mapping it to probability space, the density operator using *PBN* now reads:

In probability space: $\rho = |\Omega)P(\Omega|, \quad \mathrm{Tr}[\rho^2] = \mathrm{Tr}\rho = P(\Omega|\Omega) = 1$ (1.13b)

***Conditional density operator* (*CDO*)**: Now we are ready to define our *CDO* in both Hilbert space and probability space:

$$\rho_A = \rho \mathbf{1}_A, \quad 0 \le \mathrm{Tr}\rho_A = \mathrm{Tr}[\rho \mathbf{1}_A] \le 1 \tag{1.14}$$

It doesn't have unit trace. But, as we will see, our *CDO* will play a crucial role in our unified expressions.

***The mapping of momentum operator***: Suppose the observable is the continuous position operator, we proposed in Ref. [5] that the momentum operator $\hat{p}$ in Hilbert space is mapped to an anti-Hermitian operator $\hat{\kappa}$ in probability space, representing imaginary wave number, as follows:

$$\langle x|\hat{p}|x'\rangle = \frac{\hbar}{i}\frac{\partial}{\partial x}\delta(x-x') \to P(x|\hat{p}|x') = P(x|i\hbar\hat{\kappa}|x') = \frac{\hbar}{i}\frac{\partial}{\partial x}\delta(x-x') \tag{1.15}$$

We will use this representation when discuss non-commutative observables later (§7).

## 2. System with one discrete observable

In this section, we discuss *CE*, *CP* and *AP* of a stable probability space with one discrete observable. We assume the observable is the stationary Hamiltonian of a bounded particle (like a harmonic oscillator), which has the following discrete eigenvalues and eigenvectors (a v-basis) in the Hilbert space:

$$\text{Discrete } \varepsilon\text{-spectrum}: \quad \hat{H}|\varepsilon_i\rangle = \varepsilon_i|\varepsilon_i\rangle, \quad \langle\varepsilon_i|\varepsilon_j\rangle = \delta_{ij}, \quad \hat{I} = \sum_i |\varepsilon_i\rangle\langle\varepsilon_i| \tag{2.1}$$

The time-independent system state ket in the v-basis is given by:

$$|\Psi\rangle = \hat{I}|\Psi\rangle = \sum_i |\varepsilon_i\rangle\langle\varepsilon_i|\Psi\rangle \equiv \sum_i |\varepsilon_i\rangle c_i \tag{2.2}$$

The expectation value of $\hat{H}$ is given by (1.4a) [1-2]:

$$E(\hat{H}) \equiv \langle\hat{H}\rangle \equiv \bar{H} \equiv \langle\Psi|\hat{H}|\Psi\rangle = \sum_i \langle\Psi|\hat{H}|\varepsilon_i\rangle\langle\varepsilon_i|\Psi\rangle = \sum_i \varepsilon_i |c_i|^2 = \sum_i \varepsilon_i P(\varepsilon_i) \tag{2.3}$$

Using *PBN*, the v-basis (2.1) is mapped to *P*-basis of the probability space $(\Omega, \hat{H}, P)$:

$$\hat{H}\,|\,\varepsilon_i) = \varepsilon_i\,|\,\varepsilon_i), \quad P(\varepsilon_i\,|\,\varepsilon_j) = \delta_{ij}, \quad \hat{I} = \sum_i |\varepsilon_i)P(\varepsilon_i\,| \tag{2.4}$$

The expectation value (2.3) now can be written as:

$$\begin{aligned} E[\hat{H}] = \langle \hat{H} \rangle \equiv P(\Omega\,|\,\hat{H}\,|\,\Omega) = \sum_{\varepsilon_i} P(\Omega\,|\,\hat{H}\,|\,\varepsilon_i)P(\varepsilon_i\,|\Omega) \; = \sum_{\varepsilon_i\in\Omega} \varepsilon_i\, P(\varepsilon_i\,|\Omega) \\ = \sum_{\varepsilon_i\in\Omega} \varepsilon_i\, P(\varepsilon_i) \underset{(2.3)}{=} \sum_{\varepsilon_i\in\Omega} \varepsilon_i\,|\,c_i\,|^2 \underset{(2.3)}{=} \sum_i \langle \Psi\,|\,\hat{H}\,|\,\varepsilon_i\rangle\langle \varepsilon_i\,|\Psi\rangle \end{aligned} \tag{2.5}$$

Let $A$ be a subset of possible outcomes in the probability space:

$$A = \{\varepsilon_a, ..., \varepsilon_b\} \subset \Omega \tag{2.6}$$

From (1.5), we can write the conditional expectation value of $\hat{H}$ given $A$ as:

$$E[\hat{H}\,|\,A] \equiv P(\Omega\,|\,\hat{H}\,|\,A) = \sum_{\varepsilon_i} P(\Omega\,|\,\hat{H}\,|\,\varepsilon_i)P(\varepsilon_i\,|A) \; = \sum_{\varepsilon_i} \varepsilon_i P(\varepsilon_i\,|A) \tag{2.7}$$

Here, the conditional probability $P(\varepsilon_i\,|\,A)$ can be expressed by definition as:

$$\begin{aligned} P(\varepsilon_i\,|\,A) = \frac{P(\varepsilon_i \cap A\,|\,\Omega)}{P(A\,|\,\Omega)} \underset{\varepsilon_i\in A\subset\Omega}{=} \frac{P(\varepsilon_i\,|\,\Omega)}{P(A\,|\,\Omega)} = \frac{P(i)}{\sum_{\varepsilon_k\in\Omega} P(A\,|\,\varepsilon_k)P(\varepsilon_k\,|\,\Omega)} \\ = \frac{P(i)}{\sum_{\varepsilon_k\in A} P(\varepsilon_k\,|\,\Omega)} = \frac{|\,c_i\,|^2}{\sum_{\varepsilon_k\in A} |\,c_k\,|^2}, \quad \varepsilon_i \in A \subset \Omega \end{aligned} \tag{2.8}$$

Therefore, the *CE* in (2.7) becomes:

$$P(\Omega\,|\,\hat{H}\,|\,A) = \sum_{\varepsilon_k\in A} \varepsilon_i P(\varepsilon_i\,|\,A) = \frac{\sum_{\varepsilon_i\in A} \varepsilon_i\,|\,c_i\,|^2}{\sum_{\varepsilon_k\in A} |\,c_k\,|^2} \tag{2.9}$$

In terms of QM, Eq. (2.9) is the expected energy of the particle when the observed energy of the particle in is the range given by Eq. (2.6).

The indicator operator for set $A$ in the probability space is given by (1.9a):

$$\boldsymbol{1}_A = \sum_{\varepsilon_i\in A} |\,\varepsilon_i)P(\varepsilon_i\,| \tag{2.10}$$

Using it, we have the following expression:

$$P(\Omega\,|\,\hat{H}\boldsymbol{1}_A\,|\,\Omega) = \sum_{\varepsilon_k\in A} P(\Omega\,|\,\hat{H}\,|\,\varepsilon_k)P(\varepsilon_k\,|\,\Omega) = \sum_{\varepsilon_k\in A} \varepsilon_i P(\varepsilon_i\,|\,\Omega)$$

$$= P(A|\Omega)\frac{\sum_{\varepsilon_k\in A}\varepsilon_i P(\varepsilon_i \cap A|\Omega)}{P(A|\Omega)} \underset{(2.8)}{=} P(A|\Omega)\, E[\hat{H}\,|\,A] \tag{2.11}$$

Here, we have actually derived Eq. (1.8) for the case of discrete spectrum:

$$E[\hat{H}\,|\,A] = P(\Omega\,|\,\hat{H}\,|\,A) = \frac{P(\Omega\,|\,\hat{H}\boldsymbol{1}_A\,|\,\Omega)}{P(A|\Omega)} \tag{2.12}$$

Using density operator in (1.13a), the expectation value (2.5) can be written as:

$$E[\hat{H}] = \sum_i \langle\Psi\,|\,\hat{H}\,|\,\varepsilon_i\rangle\langle\varepsilon_i\,|\Psi\rangle = \mathrm{Tr}\{\hat{H}\,|\Psi\rangle\langle\Psi\,|\} = \mathrm{Tr}[\hat{H}\rho] = \frac{\mathrm{Tr}[\hat{H}\rho]}{\mathrm{Tr}[\rho]} \tag{2.13}$$

The indicator operator (2.10) in Hilbert space is given by:

$$\boldsymbol{1}_A = \sum_{\varepsilon_i\in A}|\,\varepsilon_i\rangle\langle\varepsilon_i\,|,\quad \boldsymbol{1}_A\,|\,\varepsilon_i\rangle = \begin{cases}0, & \text{if } \varepsilon_i \notin A\\ |\,\varepsilon_i\rangle, & \text{if } \varepsilon_i \in A\end{cases} \tag{2.14}$$

Now we can express *CE* (2.12) in terms of density and indicator operators:

$$E[\hat{H}\,|\,A] \underset{(1.8)}{=} \frac{P(\Omega\,|\,\hat{X}\boldsymbol{1}_A\,|\,\Omega)}{P(A|\Omega)} = \frac{\sum_{\varepsilon_i\in A}\varepsilon_i\,|\,c_i\,|^2}{\sum_{\varepsilon_k\in A}|\,c_k\,|^2} = \frac{\mathrm{Tr}[\hat{H}\rho\boldsymbol{1}_A]}{\mathrm{Tr}[\rho\boldsymbol{1}_A]} \tag{2.15}$$

$$\boldsymbol{CE}\text{: } E[\hat{H}\,|\,A] = \frac{P(\Omega\,|\,\hat{X}\boldsymbol{1}_A\,|\,\Omega)}{P(A|\Omega)} = \frac{\mathrm{Tr}[\hat{H}\rho_A]}{\mathrm{Tr}[\rho_A]} \tag{2.16}$$

The denominator of Eq. (2.16) represents the absolute probability (*AP*):

$$\boldsymbol{AP}\text{: } P(A) \equiv P(A|\Omega) = \sum_{\varepsilon_k\in A} P(\varepsilon_i\,|\,A) = \sum_{\varepsilon_k\in A}|\,c_k\,|^2 = \mathrm{Tr}[\rho\boldsymbol{1}_A] \underset{(1.14)}{=} \mathrm{Tr}[\rho_A] \tag{2.17}$$

By definition, the conditional probability (*CP*) now can be rewritten as:

$$\boldsymbol{CP}\text{: } P(A\,|\,B) \equiv \frac{P(A\cap B\,|\,\Omega)}{P(B)} = \underset{(2.15b)}{=} \frac{\sum_{\varepsilon_i\in A\cap B}|\,c_i\,|^2}{\sum_{\varepsilon_k\in B}|\,c_k\,|^2} =$$

$$= \frac{\sum_{\varepsilon_i\in A\cap B}|\langle\varepsilon_i\,|\Psi\rangle|^2}{\sum_{\varepsilon_k\in B}|\langle\varepsilon_k\,|\Psi\rangle|^2} = \frac{\mathrm{Tr}[\rho\boldsymbol{1}_{A\cap B}]}{\mathrm{Tr}[\rho\boldsymbol{1}_B]} = \frac{\mathrm{Tr}[\rho_{A\cap B}]}{\mathrm{Tr}[\rho_B]} \tag{2.18}$$

We will see the unified format of (2.16-18) in our next sections.

## 3. System with one continuous observable

In this section, we investigate *CE*, *CP* and *AP* of stable probability space with one continuous observable. We assume the observable is the position $\hat{X}$ of a particle, which has the following continuous spectrum (a v-basis) in the Hilbert space as in (1.2b-1.3b):

$$\text{Continuous } x\text{-spectrum}: \quad \hat{X}\,|\,x\rangle = x\,|\,x\rangle, \quad \langle x\,|\,x'\rangle = \delta(x-x'), \quad \hat{I} = \int dx\,|x\rangle\langle x| \tag{3.1}$$

The time-independent system state ket in the v-basis is given by:

$$|\,\Psi\rangle = \hat{I}\,|\,\Psi\rangle = \int dx\,|x\rangle\langle x\,|\,\Psi\rangle \equiv \int dx\,|x\rangle\,\Psi(x) \tag{3.2}$$

The expectation value of $\hat{X}$ is given by (1.4b):

$$E(\hat{X}) \equiv \langle\Psi\,|\,\hat{X}\,|\,\Psi\rangle = \int dx\,\langle\Psi\,|\,\hat{X}\,|\,x\rangle\langle x\,|\Psi\rangle = \int dx\,x\,|\,\Psi(x)\,|^2 = \int dx\,x\,P(x) \tag{3.3}$$

Using *PBN*, the v-basis (2.1) is mapped to *P*-basis in a probability space $(\Omega,\ \hat{X}, P)$, as in (1.2b-1.3b):

$$\hat{X}\,|\,x) = x\,|\,x), \quad P(x\,|\,x') = \delta(x-x'), \quad \hat{I} = \int dx\,|x)P(x| \tag{3.4}$$

The expectation value (3.3) now can be written as:

$$\begin{aligned} E(\hat{X}) &\equiv \langle\hat{X}\rangle \equiv \bar{X} = P(\Omega\,|\,\hat{X}\,|\,\Omega) = \int_{x\in\Omega} P(\Omega\,|\,\hat{X}\,|\,x)dx\,P(x\,|\,\Omega) \\ &= \int_{x\in\Omega} dx\,x\,P(x) = \int_{x\in\Omega} dx\,x\,|\,\Psi(x)\,|^2 \end{aligned} \tag{3.5}$$

Now let us define a subset of possible outcomes in the probability space:

$$A = [x_a, x_b] \subset \Omega \tag{3.6}$$

Then, from (1.6), the conditional expectation value of $\hat{H}$ given *A* is:

$$E[\hat{X}\,|\,A] \equiv P(\Omega\,|\,\hat{X}\,|\,A) = \int_{x\in\Omega} P(\Omega\,|\,\hat{X}\,|\,x)dx\,P(x\,|\,A) = \int_{x\in\Omega} dx\,x\,P(x\,|\,A) \tag{3.7}$$

Here, the conditional probability $P(x\,|\,A)$ can be expressed by definition as:

$$P(x\,|\,A) = \frac{P(x\cap A\,|\,\Omega)}{P(A\,|\,\Omega)} \underset{x\in A\subset\Omega}{=} \frac{P(x\,|\,\Omega)}{P(A\,|\,\Omega)} = \frac{P(x)}{\int_{x'\in\Omega} dxP(A\,|\,x')P(x'\,|\,\Omega)}$$

$$= \frac{P(x)}{\int_{x\in A} dx' P(x'|\Omega)} = \frac{|\Psi(x)|^2}{\int_{x\in A} dx' |\Psi(x')|^2}, \quad x \in A \subset \Omega \tag{3.8}$$

Therefore, the *CE* in (3.7) reads:

$$E[\hat{H} | A] = P(\Omega | \hat{H} | A) = \int_{x\in A} dx\, x\, P(x | A) \underset{(3.8)}{=} \frac{\int_{x\in A} dx\, x\, |\Psi(x)|^2}{\int_{x\in A} dx' |\Psi(x')|^2} \tag{3.9}$$

In terms of QM, Eq. (2.9) is the expected position of the particle when the observed position of the particle in is the range given by Eq. (3.6).

We now use the indicator operator (1.9b) in probability space $(\Omega, \hat{X}, P)$ :

$$\boldsymbol{1}_A = \int_{x\in A} dx\, |x)P(x| \tag{3.10}$$

We have following expression:

$$P(\Omega | \hat{X}\boldsymbol{1}_A | \Omega) = \int_{x\in A} dx P(\Omega | \hat{X} | x) P(x | \Omega) = \int_{x\in A} dx\, x\, P(x | \Omega)$$
$$= P(A|\Omega) \frac{\int_{x\in A} dx\, x\, P(x \cap A|\Omega)}{P(A|\Omega)} \underset{(3.8)}{=} P(A|\Omega)\, E(\hat{X} | A) \tag{3.11}$$

Here, we have actually derived Eq. (1.8) for the case of continuous spectrum:

$$E[\hat{X} | A] = P(\Omega | \hat{X} | A) = \frac{P(\Omega | \hat{X}\boldsymbol{1}_A | \Omega)}{P(A|\Omega)} \tag{3.12}$$

Using density operator in (1.13a), the expectation value in (3.3) can be written as:

$$E(\hat{X}) \equiv \int_x dx\, x\, |\Psi(x)|^2 = \int_x dx \langle \Psi | \hat{X} | x \rangle \langle x | \Psi \rangle = \mathrm{Tr}[\hat{H} | \Psi \rangle \langle \Psi |] = \mathrm{Tr}[\hat{H}\rho] \tag{3.13}$$

The indicator operator (3.10) in the Hilbert space is given by (1.9b):

$$\boldsymbol{1}_A = \int_{x\in A} dx |x\rangle\langle x|, \quad \boldsymbol{1}_A |x\rangle = \begin{cases} 0, & \text{if } x \notin A \\ |x\rangle, & \text{if } x \in \bar{a} \subset A \ \& \int_{\bar{a}} dx > 0 \\ \delta(x - x')|x\rangle, & \text{if } A = \{|x'\rangle\langle x'| \cup B\} \ \& \ x \notin B \end{cases} \tag{3.14}$$

Now we express (3.9) in terms of density and indicator operators:

$$E[\hat{X} \mid A] \underset{(1.8)}{=} \frac{P(\Omega \mid \hat{X}\boldsymbol{1}_A \mid \Omega)}{P(A \mid \Omega)} = \frac{\int_{x \in A} dx\, x \mid \Psi(x) \mid^2}{\int_{x \in A} dx \mid \Psi(x) \mid^2} = \frac{\mathrm{Tr}[\hat{X} \rho \boldsymbol{1}_A]}{\mathrm{Tr}[\rho \boldsymbol{1}_A]} \tag{3.15}$$

Once again, we have expressions with the same unified format as in (2.16-18):

$$\boldsymbol{CE}\text{:}\ E[\hat{X} \mid A] = \frac{P(\Omega \mid \hat{X}\boldsymbol{1}_A \mid \Omega)}{P(A \mid \Omega)} = \frac{\mathrm{Tr}[\hat{X} \rho_A]}{\mathrm{Tr}[\rho_A]} \tag{3.16}$$

$$\boldsymbol{AP}\text{:}\ P(A) = P(A \mid \Omega) = \int_{x \in A} dx\, P(x \mid A) = \int_{x \in A} dx \mid \Psi(x) \mid^2 = \mathrm{Tr}[\rho \boldsymbol{1}_A] = \mathrm{Tr}[\rho_A] \tag{3.17}$$

$$\boldsymbol{CP}\text{:}\ P(A \mid B) \equiv \frac{P(A \cap B \mid \Omega)}{P(B \mid \Omega)} \underset{(3.15b)}{=} \frac{\int_{x \in A \cap B} dx \mid \langle x \mid \Psi \rangle \mid^2}{\int_{x \in' B} dx' \mid \langle x' \mid \Psi \rangle \mid^2} = \frac{\mathrm{Tr}[\rho \boldsymbol{1}_{A \cap B}]}{\mathrm{Tr}[\rho \boldsymbol{1}_B]} = \frac{\mathrm{Tr}[\rho_{A \cap B}]}{\mathrm{Tr}[\rho_B]} \tag{3.18}$$

## 4. System with two commutative observables

Assume we have a 2D-particle, its positions are $X$ and Y respectively and a stable probability space $(\Omega, \Im, P)$ is associated with them [5]. Then there exists a joint stable probability density $P(x, y)$ for any joint event $(x \cap y) \in \Im$, such that:

$$| x, y) \equiv | x \cap y) \equiv | X = x \cap Y = y), \quad X | x, y) = x | x, y), \quad Y | x, y) = y | x, y) \tag{4.1a}$$

$$P(\Omega | x, y) = 1, \quad | \Omega) \equiv | \Omega_{x,y}), \quad P(\Omega | = P(\Omega_{x,y} |, \tag{4.1b}$$

$$P(x, y | x', y') = \delta(x - x')\delta(y - y'), \quad P(x, y | \Omega) \equiv P(x, y) \geq 0 \tag{4.1c}$$

The stable probability density $P(x, y)$ has following additional properties:

$$\int | x, y)\, dxdy P(x, y | = \hat{I}, \quad P(\Omega | \Omega) = \int dxdy\, P(x, y | \Omega) = \int dxdy\, P(x, y) = 1 \tag{4.2}$$

$$P(x | \Omega) \equiv P(x, * | \Omega) = \int dy\, P(x, y | \Omega), \quad P(y | \Omega) \equiv P(y, * | \Omega) = \int dx\, P(x, y | \Omega) \tag{4.3}$$

Because the basis is complete, if $A$ is any set in $\Im$, we have:

$$\begin{aligned} XY | A) &= \int XY | x, y)\, dxdy\, P(x, y | A) = \int xy | x, y) dxdy P(x, y | A) = YX | A) \\ \therefore (XY - YX) | A) &\equiv [X, Y] | A) = 0, \quad \textit{for any } A \in \Im \end{aligned} \tag{4.4}$$

Here we have used the definition of *commutator* in QM (see [6], §7.3). As we proposed in Ref [3], the existence of joint density means that, in the induced Hilbert space, the corresponding operators must have a complete common eigenvectors, and therefore they must be *commutative* operators:

$$\hat{X}\,|\,x,y\rangle = x\,|\,x,y\rangle,\;\; \hat{Y}\,|\,x,y\rangle = y\,|\,x,y\rangle \\ [\hat{X}-\hat{Y}] \equiv \hat{X}\hat{Y}-\hat{Y}\hat{X}=0 \tag{4.5}$$

The v-basis in the induced Hilbert space of the 2D QM problem is mapped from (4.1-2):

$$|\,x,y\rangle \equiv |\,x\cap y\rangle \equiv |\,X=x\cap Y=y\rangle, \quad X\,|\,x,y\rangle = x\,|\,x,y\rangle, \quad Y\,|\,x,y\rangle = y\,|\,x,y\rangle \\ \langle x,y\,|\,x',y'\rangle = \delta(x-x')\delta(y-y'), \quad P(x,y\,|\,\Omega) = P(x,y) \equiv |\,\Psi(x,y)\,|^2 \geq 0 \tag{4.6}$$

$$\int |\,x,y\rangle\, dxdy\, \langle x,y\,| = \hat{I}, \quad \langle \Psi\,|\,\Psi\rangle = \int dxdy\, P(x,y) = \int dxdy\,|\,\Psi(x,y)\,|^2 = 1 \tag{4.7}$$

$$P(x) \equiv P(x,*) = \int dy\, P(x,y) = \int dy\,|\,\Psi(x,y)\,|^2 \\ P(y) \equiv P(*,y) = \int dx\, P(x,y) = \int dx\,|\,\Psi(x,y)\,|^2 \tag{4.8}$$

Therefore, (4.3) can be expressed using probability density in Hilbert space:

$$P(x\,|\,\Omega) \equiv P(x,*\,|\,\Omega) = \int dy\, P(x,y\,|\,\Omega) = \int dy\, P(x,y) = \int dy\,|\,\Psi(x,y)\,|^2, \\ P(y\,|\,\Omega) \equiv P(*,y\,|\,\Omega) = \int dx\, P(x,y\,|\,\Omega) = \int dy\, P(x,y) = \int dy\,|\,\Psi(x,y)\,|^2 \tag{4.9}$$

Assume *g*(*x*, *y*) is a Borel function [8]. For $\forall A \in \Im$, we have the following conditional probability of g(*X*, *Y*) given *A*:

$$\boldsymbol{CE}\text{: } E[g(\hat{X},\hat{Y})\,|\,A] = P(\Omega\,|\,g(\hat{X},\hat{Y})\,|\,A) = \int_{x,y\in H} dx\,dy\, g(x,y)\, P(x,y\,|\,A) = \\ \frac{\int_{x,y\in A} dx\,dy\, g(x,y)\, P(x,y\,|\,\Omega)}{\int_{x,y\in A} dx\,dy\, P(x,y\,|\,\Omega)} = \frac{\int_{x,y\in A} dx\,dy\, g(x,y)\,|\,\Psi(x,y)\,|^2}{\int_{x,y\in A} dx\,dy\,|\,\Psi(x,y)\,|^2} \\ = \frac{\mathrm{Tr}\left\{g(\hat{X},\hat{Y})\rho \boldsymbol{1}_A\right\}}{\mathrm{Tr}\left\{\rho \boldsymbol{1}_A\right\}} = \frac{\mathrm{Tr}\left\{g(\hat{X},\hat{Y})\rho_A\right\}}{\mathrm{Tr}\left\{\rho_A\right\}} \tag{4.10}$$

We also have following expressions for *AP* and *CP*:

$$\boldsymbol{AP}\text{: } P(A) = P(A\,|\,\Omega) = \int_{x,y\in A} dx\,dy P(x,y\,|\,\Omega) = \int_{x,y\in A} dx\,dy\,|\,\Psi(x,y)\,|^2 = \mathrm{Tr}\left\{\rho \boldsymbol{1}_A\right\} \tag{4.11}$$

$$\boldsymbol{CP}\text{: } P(A\,|\,B) \equiv \frac{P(A\cap B\,|\,\Omega)}{P(B\,|\,\Omega)} = \frac{\int_{(x,y)\in A\cap B} dx\,|\,\langle x,y\,|\,\Psi\rangle\,|^2}{\int_{(x',y')\in B} dx'\,|\,\langle x',y'\,|\,\Psi\rangle\,|^2} = \frac{\mathrm{Tr}[\rho \boldsymbol{1}_{A\cap B}]}{\mathrm{Tr}[\rho \boldsymbol{1}_B]} = \frac{\mathrm{Tr}[\rho_{A\cap B}]}{\mathrm{Tr}[\rho_B]}. \tag{4.12}$$

We see that (4.10-12) have the same unified format as in (2.16-18) and (3.16-18).

***Two independent observables***: Now let us assume that *X* and *Y* are the positions of 2-dimensional particle with separable potential (see Eq. (4.17) below). Then *X* and *Y* are *two independent observables* [4] on the probability space $(\Omega, \Im, P)$ and there exist two subsets $\Omega_x$ and $\Omega_y$ such that $\Omega = \Omega_x \otimes \Omega_y$, $\Omega_x \cap \Omega_y = \varnothing$, or:

$$|\Omega) = |\Omega_{x,y}) = |\Omega_x)\,|\Omega_y), \quad P(\Omega| = P(\Omega_{x,y}| = P(\Omega_x|\,P(\Omega_y| \tag{4.13}$$

For any joint event $(x, y) \equiv (X = x \wedge Y = y) \equiv (X = x \otimes Y = y) \equiv (x \otimes y) \in \Im$, we have:

$$|x, y) = |x)\,|y), \quad P(x, y| = P(x|\,P(y| \tag{4.14}$$

$$P(x|\Omega_x) \equiv f_X(x), \quad P(y|\Omega_y) \equiv f_Y(x) \tag{4.15a}$$

$$P(\Omega|x, y) = 1, \quad P(x, y|\Omega) = P(x|\Omega_x)P(y|\Omega_y) = f_X(x)\,f_Y(y) \tag{4.15b}$$

Therefore, for $\forall x \in \Omega_x$, we have following *CE* of a *Borel* function *g*(*X*) given *X* = *x*:

$$P(\Omega|\,g(X)\,|x) = \int P(\Omega|\,g(X)\,|x')\,dx'\,P(x'|x) = \int g(x')dx'\,\delta(x - x') = g(x) \tag{4.16}$$

In Hilbert space, the quantum system has following 2D-Hamiltonian:

$$\hat{H} = \hat{H}_x + \hat{H}_y = \frac{\hat{p}_x^{\,2}}{2m} + V_1(x) + \frac{\hat{p}_y^{\,2}}{2m} + V_2(y) \tag{4.17}$$

Hence the system wave function is the product of wave function of each dimension:

$$\langle x, y|\Psi\rangle = \langle x|\Psi_1\rangle\langle y|\Psi_2\rangle = \Psi_1(x)\Psi_2(y) \tag{4.18}$$

And the density operator can also be written as a product:

$$\rho = |\Psi\rangle\langle\Psi| = \{|\Psi_1\rangle\langle\Psi_1|\}\{|\Psi_2\rangle\langle\Psi_2|\} = \rho_1\,\rho_2 \tag{4.19}$$

This leads to the following absolute probabilities:

$$P(x, y) \equiv |\Psi(x, y)|^2 = |\Psi_1(x)|^2\,|\Psi_2(y)|^2 \tag{4.20}$$

$$\begin{aligned} P(x) &\equiv P(x,*) = \int dy\,P(x, y) = \int dy\,|\Psi(x, y)|^2 = |\Psi_1(x)|^2 \\ P(y) &\equiv P(*, y) = \int dx\,P(x, y) = \int dx\,|\Psi(x, y)|^2 = |\Psi_2(y)|^2 \end{aligned} \tag{4.21}$$

Because the subset *A* is also separable, $|A) = |A_1)\,|A_2)$, we get *CE* in our general format:

$$E[\hat{X} \mid A] = \frac{P(\Omega \mid \hat{X}\boldsymbol{1}_A \mid \Omega)}{P(H \mid \Omega)} = \frac{\int_{x \in A_1} dx\ x \mid \Psi_1(x) \mid^2}{\int_{x \in A_1} dx \mid \Psi_1(x) \mid^2} = \frac{tr\left\{\hat{X}\rho_1 \boldsymbol{1}_{A_1}\right\}}{tr\left\{\rho_1 \boldsymbol{1}_{A_1}\right\}} = \frac{tr\left\{\hat{X}\rho_{A_1}\right\}}{tr\left\{\rho_{A_1}\right\}} \tag{4.22a}$$

$$E[\hat{Y} \mid A] = \frac{P(\Omega \mid \hat{X}\boldsymbol{1}_A \mid \Omega)}{P(H \mid \Omega)} = \frac{\int_{y \in A_2} dy\ y \mid \Psi_2(y) \mid^2}{\int_{y \in A_2} dy \mid \Psi_2(y) \mid^2} = \frac{tr\left\{\hat{Y}\rho_2 \boldsymbol{1}_{A_2}\right\}}{tr\left\{\rho_2 \boldsymbol{1}_{A_2}\right\}} = \frac{tr\left\{\hat{Y}\rho_{A_2}\right\}}{tr\left\{\rho_{A_2}\right\}} \tag{4.22b}$$

They also share our unified *CE* format in (2.16), (3.16) and (4.10).

***Two commutable but dependent observables:*** If *X* and *Y* are commutative but *not independent* of each other, we cannot make the assertions (4.12-22). In this case, to calculate *CDO* in (4.10), we have to use $P(x, y) = \mid \Psi(x, y) \mid^2$, that is:

$$\begin{aligned} tr\left\{g(\hat{X},\hat{Y})\rho_A\right\} &= \int_{(x',y')\in A} dx\,dy\,dx'\,dy' \langle x, y \mid g(\hat{X},\hat{Y}) \mid \Psi\rangle\langle\Psi \mid x', y'\rangle \langle x', y' \mid x, y\rangle \\ &= \int_{(x',y')\in A} dx\,dy\,dx'\,dy'\, g(x,y)\langle x, y \mid \Psi\rangle\langle\Psi \mid x', y'\rangle\, \delta(x'-x)\delta(y'-y) \\ &= \int_{(x,y)\in A} dx\,dy\, g(x,y) \mid \Psi(x,y) \mid^2 \end{aligned} \tag{4.23}$$

Here $\Psi(x, y)$ is the exact solution of the following 2D-Schrodinger equation:

$$\hat{H} = \hat{H}_x + \hat{H}_y + V(x, y) = \frac{\hat{p}_x^{\ 2}}{2m} + V_1(x) + \frac{\hat{p}_y^{\ 2}}{2m} + V_2(y) + V(x, y) \tag{4.24}$$

## 5. Many particle system and Fock space

In our previous work (see §3.1 of Ref [3]), we discussed the expectation values for Fock space (§22, [6]) of system of indistinguishable non-interacting particles using *PBN*). The basis vectors are eigenvectors of occupation number operators $\hat{\vec{N}} = (\hat{n}_1, \ldots, \hat{n}_t)$:

$$[\hat{n}_i, \hat{n}_k] = 0, \quad \hat{n}_i \left|\vec{N}\right\rangle = \hat{n}_i \mid n_1, n_2, \ldots, n_t\rangle = n_i \mid n_1, n_2, \ldots, n_t\rangle \tag{5.1}$$

$$\left|\vec{N}\right\rangle = \mid n_1, n_2, \ldots, n_t\rangle = \prod_{k=1}^{t} \mid n_k\rangle, \quad \mid n_j\rangle \mid n_k\rangle = \mid n_k\rangle \mid n_j\rangle, \quad [\hat{n}_i, \hat{n}_j] = 0 \tag{5.2}$$

They are mapped to the basis of probability space:

$$\mid \vec{N}) = \mid n_1, n_2, \ldots, n_t) = \prod_{i=1}^{t} \mid n_i)_i, \quad \mid n_i)_i \mid n_j)_j = \mid n_j)_j \mid n_i)_i \tag{5.3}$$

Theses *P*-kets are eigenvectors of the random observables $\hat{\boldsymbol{N}} = (\hat{N}_1, \ldots, \hat{N}_t)$:

$$\hat{N}_i \,|\, n_1,\ldots,n_t) = n_i \,|\, n_1,\ldots,n_t) \tag{5.4}$$

The two spaces share observables $\hat{\boldsymbol{N}}$ and have equivalent properties:

$$\langle \vec{N}' \,|\, \vec{N}\rangle = \delta_{\vec{N}',\vec{N}}, \sum_{\vec{N}} |\, \vec{N}\rangle\langle \vec{N} \,| = \hat{I} \Leftrightarrow P(\vec{N}' \,|\, \vec{N}) = \delta_{\vec{N}',\vec{N}}, \sum_{\vec{N}} |\, \vec{N})P(\vec{N} \,| = \hat{I} \tag{5.5}$$

Now let us study many-particle systems in Thermophysics. From quantum statistics (see [10], §4 and §5), we know that the grand partition function of a system of many identical particles is defined in terms of total energy $E_j$ and total occupation number $N$ as:

$$Z_G = \sum_{N,j} \exp[-\beta(E_j - \mu N)] = \sum_{N,j} \langle N, j \,|\, \exp[-\beta(\hat{H} - \mu\hat{N})] \,|\, N, j\rangle$$
$$= \mathrm{Tr}(\exp[-\beta(\hat{H} - \mu\hat{N})] \tag{5.6}$$

For any operator $\hat{O}$, the ensample average $\langle \hat{O} \rangle$ is obtained by the representation:

$$\langle \hat{O}\rangle = \frac{\mathrm{Tr}\{\exp[-\beta(\hat{H} - \mu\hat{N})]\hat{O}\}}{\mathrm{Tr}\{\exp[-\beta(\hat{H} - \mu\hat{N})]\}} \tag{5.7}$$

In Fock space, the total Hamiltonian and the operator of total occupation number are:

$$\hat{H} = \sum \hat{n}_j \varepsilon_j, \quad \hat{N} = \sum_j \hat{n}_j, \quad \hat{n}_j \,|\, \vec{N}\rangle = \hat{n}_j \,|\, n_1,\ldots n_\infty\rangle = n_j \,|\, n_1,\ldots n_\infty\rangle \tag{5.8}$$

Using Eq. (5.2), we can factor the grand partition function as (see [10], page 37):

$$Z_G = \sum_{\vec{N}} \langle \vec{N} \,|\, \exp[-\beta(\hat{H} - \mu\hat{N})] \,|\, \vec{N}\rangle = \prod_{i=1}^{\infty} \sum_{n_i} \langle n_i \,| \exp[-\beta(\varepsilon_i - \mu)n_i]\,|\, n_i\rangle$$
$$= \prod_{i=1}^{\infty} \mathrm{Tr}\{\exp[-\beta(\varepsilon_i - \mu)\hat{n}_i]\} = \prod_{i=1}^{\infty} Z_i \tag{5.9}$$

In Fock space, the total Hamiltonian and the operator of total occupation number are:

$$\hat{H} = \sum \hat{n}_j \varepsilon_j, \quad \hat{N} = \sum_j \hat{n}_j, \quad \hat{n}_j \,|\, \vec{N}\rangle = \hat{n}_j \,|\, n_1,\ldots n_\infty\rangle = n_j \,|\, n_1,\ldots n_\infty\rangle \tag{5.10}$$

If an operator is a linear function of occupation numbers in the following form:

$$O(\hat{\vec{N}}) = \sum_{i=1}^{\infty} a_i \hat{n}_i \tag{5.11}$$

Then its expectation value can be obtained as:

$$\langle O(\hat{\vec{N}})\rangle = \sum_{i=1}^{\infty} \frac{\mathrm{Tr}\{a_i \hat{n}_i \exp[-\beta(\varepsilon_i - \mu)\hat{n}_i]\}}{Z_i} \tag{5.12}$$

Using *PBN*, the probability of the system at one-particle-state *j* with energy $\varepsilon_j$ and occupation number $n_j$ is given by (see also [9], §11.6):

$$P(n_j \mid \Omega_j) = m(n_j) = \frac{\exp[-(n_j\varepsilon_j - \mu n_j)/kT]}{Z_j} \tag{5.13}$$

Now we are ready to find the equilibrium system state at temperature *T* in Fock space:

$$|\Psi\rangle = \sum_{\vec{N}} C(\vec{N}) \,|\, \vec{N}\rangle = \sum_{\vec{N}} C(n_1, n_2, \ldots) \,|\, n_1, n_2, \ldots\rangle = \sum_{\vec{N}} \prod_{j=1} c(n_j) \,|\, n_j\rangle \tag{5.14}$$

Based on Eq. (5.13), we can express the coefficients as:

$$C(\vec{N}) = \sqrt{(\vec{N} \mid \Omega)} = \sqrt{\prod_{j=1}^{\infty} (n_j \mid \Omega)} = \prod_{j=1}^{\infty} \frac{\exp[-(n_j\varepsilon_j - \mu n_j)/2kT]}{\sqrt{Z_j}} \tag{5.15}$$

This leads to the well known general format of expected value in Fock space [10]:

$$\langle \Psi \mid \hat{O} \mid \Psi \rangle = \sum_{\vec{N}} O(\vec{N}) C^2(\vec{N}) = \sum_{\vec{N}_j} O(\vec{N}) P(\vec{N}) = P(\Omega \mid \hat{O} \mid \Omega) = \mathrm{Tr}[\hat{O}\rho] \tag{5.16}$$

Now we can express *CE* of $\hat{O}$ given *A* as a subset of Ω in our general format:

$$\boldsymbol{CE}\text{: } E[\hat{O} \mid A] = \frac{P(\Omega \mid \hat{O}\boldsymbol{1}_A \mid \Omega)}{P(A \mid \Omega)} = \frac{\mathrm{Tr}\left\{\hat{O}\rho\boldsymbol{1}_A\right\}}{\mathrm{Tr}\left\{\rho\boldsymbol{1}_A\right\}} = \frac{\mathrm{Tr}\left\{\hat{O}\rho_A\right\}}{\mathrm{Tr}\left\{\rho_A\right\}} \tag{5.17}$$

$$\rho_A = \rho\boldsymbol{1}_A, \quad \boldsymbol{1}_A = \int_{\vec{N}\in A} d\vec{N} \,|\, \vec{N}\rangle\langle\vec{N} \,| \tag{5.18}$$

## 6. Time-dependent unified expressions of *CE* and *CP*

We have investigated various *stable* **(**time-independent**)** probability spaces, and derived following unified expressions for *CE*, *AP* and *CP* in terms of *CDO*:

$$\boldsymbol{CE}\text{: } E[\hat{O} \mid A] = \frac{P(\Omega \mid \hat{O}\boldsymbol{1}_A \mid \Omega)}{P(A \mid \Omega)} = \frac{\mathrm{Tr}\left\{\hat{O}\rho\boldsymbol{1}_A\right\}}{\mathrm{Tr}\left\{\rho\boldsymbol{1}_A\right\}} = \frac{\mathrm{Tr}\left\{\hat{O}\rho_A\right\}}{\mathrm{Tr}\left\{\rho_A\right\}} \tag{6.1}$$

$$\boldsymbol{AP}\text{: } P(A) = P(A\,|\,\Omega) = \mathrm{Tr}[\rho \boldsymbol{1}_A] = \mathrm{Tr}[\rho_A] \tag{6.2}$$

$$\boldsymbol{CP}\text{: } P(A\,|\,B) \equiv \frac{P(A \cap B\,|\,\Omega)}{P(B\,|\,\Omega)} = \frac{\mathrm{Tr}[\rho_{A\cap B}]}{\mathrm{Tr}[\rho_B]} = \frac{\mathrm{Tr}[\rho \boldsymbol{1}_A \boldsymbol{1}_B]}{\mathrm{Tr}[\rho \boldsymbol{1}_B]} \tag{6.3}$$

They have the same format for discrete or continuous spectrum, and can be interpreted in both Hilbert space and probability space:

In Hilbert space: $\rho_A = \rho \boldsymbol{1}_A, \quad O(\hat{\xi})\,|\,\xi\rangle = O(\xi)\,|\,\xi\rangle, \quad \rho = |\,\Psi\rangle\langle\Psi\,|$ (6.4)

In probability space: $\rho_A = \rho \boldsymbol{1}_A, \quad O(\hat{\xi})\,|\,\xi) = O(\xi)\,|\,\xi), \quad \rho = |\,\Omega)P(\Omega\,|$ (6.5)

Note the right most expression in (6.3). It may be used when sets $A$ and $B$ are defined in different bases associated with non-commutative observables (see §7).

The time-evolution of the system is determined by the Schrodinger equation [6-7]:

$$i\hbar \frac{\partial}{\partial t}\,|\,\Psi(t)\rangle = \hat{H}(t)\,|\,\Psi(t)\rangle, \quad \hat{H}^\dagger(t) = \hat{H}(t) \tag{6.6}$$

It has the following symbolic time-dependent solution:

$$|\,\Psi(t)\rangle = \hat{U}(t)\,|\,\Psi(0)\rangle, \qquad \hat{U}(t) = e^{i\int_0^t dt\,\hat{H}\,t/\hbar}, \quad \hat{U}^\dagger(t)\hat{U}(t) = 1 \tag{6.7}$$

We are not interested in solving this equation and just assume that we have the exact solution. Now our *CDO* in (6.4-5) are extended to following time-dependent formulas in Hilbert space or in probability space [1-2]:

In Hilbert space: $\rho_A(t) = \rho(t)\boldsymbol{1}_A, \quad O(\hat{\xi})\,|\,\xi\rangle = O(\xi)\,|\,\xi\rangle, \quad \rho(t) = |\,\Psi(t)\rangle\langle\Psi(t)\,|$ (6.8)

In probability space: $\rho_A(t) = \rho(t)\boldsymbol{1}_A, \quad O(\hat{\xi})\,|\,\xi) = O(\xi)\,|\,\xi), \quad \rho(t) = |\,\Omega_t)P(\Omega\,|$ (6.9)

In both spaces: $\mathrm{Tr}\rho(t) = \mathrm{Tr}[\rho(t)\hat{I}] = \mathrm{Tr}[\rho(0)] = 1$ (6.10)

It is easy to verify (6.10). In Hilbert space, we have:

$$\mathrm{Tr}\rho(t) = \mathrm{Tr}\big[|\,\Psi(t)\rangle\langle\Psi(t)\,|\big] = \langle\Psi(t)\,|\,\Psi(t)\rangle = \langle\Psi(0)\,|\,U^\dagger(t)U(t)\,|\,\Psi(0)\rangle = 1 \tag{6.11}$$

In probability space, we have (see Eq. (5.1.2a) of [1] and Eq. (3.4) of [2]):

$$\mathrm{Tr}\rho(t) = \mathrm{Tr}\big[|\,\Omega_t)\,P(\Omega\,|\big] = P(\Omega\,|\,\Omega_t) = 1 \tag{6.12}$$

It is also easy to show that, for both Hilbert space and probability space, we have

$$\mathrm{Tr}[\rho(t)^2] = \mathrm{Tr}\rho(t) = \mathrm{Tr}\rho(0) = 1 \tag{6.13}$$

Next, we extend the *PDF* relations in (1.6) to time-dependent ones (see Eq. (3.20) of [2]):

$$P(x_i \mid \Omega_t) = P(x_i, t) = |\langle x_i \mid \Psi(t)\rangle|^2 = |\Psi(x_i, t)|^2, \quad \text{discrete spectrum} \tag{6.14}$$

$$P(x \mid \Omega_t) = P(x, t) = |\langle x \mid \Psi(t)\rangle|^2 = |\Psi(x, t)|^2, \quad \text{cotinuous spectrum} \tag{6.15}$$

Now we can extend (6.1-3) to time-dependent unified expressions as follows:

$$\textbf{\textit{CE}}\text{: } E[\hat{O} \mid A](t) = \frac{P(\Omega \mid \hat{O}\boldsymbol{1}_A \mid \Omega_t)}{P(A \mid \Omega_t)} = \frac{\mathrm{Tr}\left\{\hat{O}\rho_A(t)\right\}}{\mathrm{Tr}\left\{\rho_A(t)\right\}} \tag{6.16}$$

$$\textbf{\textit{AP}}\text{: } P(A,t) = P(A \mid \Omega_t) = \mathrm{Tr}[\rho(t)\boldsymbol{1}_A] = \mathrm{Tr}[\rho_A(t)] \tag{6.17}$$

$$\textbf{\textit{CP}}\text{: } P(A \mid B)(t) \equiv \frac{P(A \cap B \mid \Omega_t)}{P(B \mid \Omega_t)} = \frac{\mathrm{Tr}[\rho(t)_{A\cap B}]}{\mathrm{Tr}[\rho(t)_B]} = \frac{\mathrm{Tr}[\rho(t)\boldsymbol{1}_A\boldsymbol{1}_B]}{\mathrm{Tr}[\rho(t)\boldsymbol{1}_B]}. \tag{6.18}$$

As an example, for the case of two commutative variables, the expressions (4.10-11) now are presented with time-dependent probability density:

$$E[g(\hat{X},\hat{Y}) \mid A](t) = \frac{\mathrm{Tr}\left\{g(\hat{X},\hat{Y})\rho_A(t)\right\}}{\mathrm{Tr}\left\{\rho_A(t)\right\}} = \frac{\int_{x,y\in A} dx\,dy\, g(x,y)\,|\Psi(x,y,t)|^2}{\int_{x,y\in A} dx\,dy\,|\Psi(x,y,t)|^2} \tag{6.12}$$

$$P(A,t) = P(A \mid \Omega_t) = \mathrm{Tr}\left\{\rho(t)_A\right\} = \int_{x,y\in A} dx\,dy\,|\Psi(x,y,t)|^2 \tag{6.13}$$

$$P(A \mid B)(t) = \frac{\mathrm{Tr}[\rho_{A\cap B}(t)]}{\mathrm{Tr}[\rho_B(t)]} = \frac{\int_{(x,y)\in A\cap B} dx\,|\Psi(x,y,t)|^2}{\int_{(x',y')\in B} dx'\,|\Psi(x',y',t)\rangle|^2} \tag{6.14}$$

## 7. Measurements of non-commutative observables

So far, we have only considered commutative observables. But, even in the simple example of 1D harmonic oscillator of QM, there are three observables (Hermitian operators): the position ( $\hat{x}$ ), the momentum ( $\hat{p}$ ), and the energy ( $\hat{H}$ ). They don't commute, so cannot be observed simultaneously. For example, from the commutator of $\hat{x}$ and $\hat{p}$ , the famous Heisenberg uncertainty relation (see §4.3 of [6]) is derived:

$$[\hat{x}, \hat{p}] = i\hbar \quad \Rightarrow \quad \Delta x\,\Delta p \geq \frac{\hbar}{2} \tag{7.1}$$

To see this in another way, let us calculate the conditional expectation of momentum given fixed position, starting from our general formula (6.1):

$$E[\hat{O} \mid A] \Rightarrow E[\hat{p} \mid x] = P(\Omega \mid \hat{p} \mid x) = \frac{\mathrm{Tr}\{\hat{p}\rho \boldsymbol{1}_x\}}{\mathrm{Tr}\{\rho \boldsymbol{1}_x\}} \tag{7.2}$$

We present our calculation of (7.2) in Hilbert space using Dirac notation as follows:

$$P(\Omega \mid \hat{p} \mid x) = \frac{\mathrm{Tr}\{\hat{p}\rho \boldsymbol{1}_x\}}{\mathrm{Tr}\{\rho \boldsymbol{1}_x\}} = \frac{\int dx'\, dx'' dx''' \langle x' \mid \hat{p} \mid x''\rangle \langle x'' \mid \Psi\rangle \langle \Psi \mid x'''\rangle \langle x''' \mid \boldsymbol{1}_x \mid x'\rangle}{\int dx' dx'' \langle x' \mid \Psi\rangle \langle \Psi \mid x''\rangle \langle x'' \mid \boldsymbol{1}_x \mid x'\rangle}$$

$$= \frac{\int dx' dx'' dx''' \langle x' \mid \hat{p} \mid x''\rangle \Psi(x'') \Psi^*(x''') \delta(x - x') \delta(x - x''')}{\int dx' dx'' \Psi(x') \Psi^*(x'') \delta(x - x') \delta(x' - x'')}$$

$$= \frac{\int dx' dx'' \frac{\hbar}{i} \frac{\partial}{\partial x'} \delta(x' - x'') \delta(x - x') \Psi(x'') \Psi^*(x)}{|\Psi(x)|^2} = \frac{\int dx' \frac{\hbar}{i} \frac{\partial}{\partial x'} \delta(x - x') |\Psi(x')|^2}{|\Psi(x)|^2}$$

$$= \int dx' \frac{\hbar}{i} \frac{\partial}{\partial x'} \delta(x - x') \quad \Rightarrow \text{undefined at } x' = x \tag{7.3}$$

We can also calculate (7.2) using *PBN* with operator $\hat{p}$ presented in (1.15):

$$P(\Omega \mid \hat{p} \mid x) = \int dx' dx'' P(\Omega \mid x') P(x' \mid \hat{p} \mid x'') P(x'' \mid x) = \int dx' dx'' \frac{\hbar}{i} \frac{\partial}{\partial x'} \delta(x' - x'') \delta(x'' - x)$$

$$= \int dx' \frac{\hbar}{i} \frac{\partial}{\partial x'} \delta(x' - x) \Rightarrow \text{undefined at } x' = x \tag{7.4}$$

Results (7.3) and (7.4) are identical and both are singular, as expected from Heisenberg uncertainty in (7.1). Therefore, our proposed mapping (1.15) and our universal expression (6.1) are self-consistent, even for singular *CE*.

We can also calculate *CP* of $p$ given $x$. We start from the right most expression of (6.3), because $p$ and $x$ belong to non-commutative operators:

$$P(p \mid x) = \frac{\mathrm{Tr}[\rho \boldsymbol{1}_x \boldsymbol{1}_p]}{\mathrm{Tr}[\rho \boldsymbol{1}_x]} = \frac{\int dx'\, dx'' \langle x' \mid \Psi\rangle \langle \Psi \mid p\rangle \langle p \mid x''\rangle \langle x'' \mid \boldsymbol{1}_x \mid x'\rangle}{\int dx' dx'' \langle x' \mid \Psi\rangle \langle \Psi \mid x''\rangle \langle x'' \mid \boldsymbol{1}_x \mid x'\rangle} = \frac{\langle \Psi \mid p\rangle \langle p \mid x\rangle}{\Psi^*(x)} \tag{7.5}$$

This is not a real quantity, so is not a *CP* by nature. But it may contain useful intermediate information, as shown by following integration:

$$\int dp\, P(p\,|\,x) = \frac{\int dp\langle\Psi\,|\,p\rangle\langle p\,|\,x\rangle}{\Psi^*(x)} = \frac{\Psi^*(x)}{\Psi^*(x)} = 1 \tag{7.6}$$

# 8. Summary and discussion

With the help of indicator operator and conditional density operator (*CDO*), we have investigated various quantum systems and have derived unified and simple expressions for *CE*, *AP* and *CP* (6.1-3), and then their time-dependent versions (6.16-17).

We call these expressions *unified*, because they have the same format for discrete or continuous spectrum and they are defined in both Hilbert space (using Dirac notation) and in probability space (using *PBN*).

We call these expressions *simple*, because they don't refer to Measure theory or Von Neumann algebra (W*-algebra) [11-12], which are beyond the general math-level of most students learning introductory quantum mechanics (QM).

However, although our expressions are simple, they can be identified with the expressions of *CDO* and *CP* in current literature. For example, if *X* and *Y* are *commutative*, the *CDO* defined in Eqs. (2.9-2.11) of Ref. [13] is

$$\rho\,|_Y \equiv \frac{\rho E(Y)}{\mathrm{Tr}\{\rho E(Y)\}} \underset{\substack{our \\ notation}}{\equiv} \frac{\rho \boldsymbol{1}_Y}{\mathrm{Tr}\{\rho \boldsymbol{1}_Y\}} = \frac{\rho_Y}{\mathrm{Tr}\{\rho_Y\}} \tag{8.1}$$

It has unit trace and is used to calculate the *CP* of *X* given *Y* as follows:

$$(\forall Y) \quad P(X\,|\,Y) = Tr[\rho_Y E(X)] = \frac{Tr[\rho E(x)E(y)]}{Tr[\rho E(y)]} \underset{\substack{our \\ notation}}{=} \frac{\mathrm{Tr}[\rho \boldsymbol{1}_X \boldsymbol{1}_Y]}{\mathrm{Tr}[\rho \boldsymbol{1}_Y]} \tag{8.2}$$

We see that, after transformed to our notation, Eq. (8.2) is identical to our Eq. (6.3). And we could have defined our *CDO* as in Eq. (8.1). But this might not be useful when expressing *AP* (absolute probability), as in (6.2) or (6.13).

Now we demonstrate how to evaluate (6.1) using both Dirac notation and *PBN.* Assuming *X* and *Y* are commutative, we apply Eq. (6.1) to the *CE* of *X* given *Y*:

$$E[X\,|\,Y] = P(\Omega\,|\,X\,|\,Y) = \frac{\mathrm{Tr}\{\hat{X}\rho_Y\}}{\mathrm{Tr}\{\rho_Y\}} \tag{8.3}$$

First, we evaluate the numerator of (8.3) in Hilbert space, using Dirac notation with the v-basis give in (4.5-7):

$$\mathrm{Tr}\left\{\hat{X}\rho_Y\right\} = \int_{x'\in\Omega} dx\,dy\,dx'\langle x,y\,|\,\hat{X}\,|\,\Psi\rangle\langle\Psi\,|\,x',Y\rangle\langle x',Y\,|\,x,y\rangle$$
$$= \int_{x'\in\Omega} dx\,dy\,dx'\,x\langle x,y\,|\,\Psi\rangle\langle\Psi\,|\,x',Y\rangle\,\delta(x'-x)\,\delta(y-Y) == \int_{x\in\Omega} dx\;\, x|\,\Psi(x,Y)\,|^2 \quad (8.4)$$

Then we evaluate the denominator of Eq. (8.3) in probability space, using *PBN* with the *P*-basis given in (4.1-2):

$$\mathrm{Tr}\left\{\rho_Y\right\} = \int_{x'\in\Omega} dx\,dy\,dx'(x,y\,|\,\Omega)P(\Omega\,|\,x',Y)P(x',Y\,|\,x,y)$$
$$= \int_{x'\in\Omega} dx\,dy\,dx'\,x\,P(x,y\,|\,\Omega)\,\delta(x'-x)\,\delta(y-Y) \underset{(4.6)}{=} \int_{x\in\Omega} dx\;\, x|\,\Psi(x,Y)\,|^2 \quad (8.5)$$

Inserting (8.4-5) into (8.3), we obtain the *CE* expression:

$$\boldsymbol{CE}\text{: } E[X\,|\,Y] = P(\Omega\,|\,X\,|\,Y) = \frac{\mathrm{Tr}\left\{\hat{X}\rho_Y\right\}}{\mathrm{Tr}\left\{\rho_Y\right\}} = \frac{\int_{x\in\Omega} dx\,x|\,\Psi(x,Y)\,|^2}{\int_{x\in\Omega} dx\,|\,\Psi(x,Y)\,|^2} \quad (8.6)$$

The *AP* defined in (6.2) is already given by (8.5):

$$\boldsymbol{AP}\text{: } P(Y) = P(Y\,|\Omega) = \mathrm{Tr}[\rho_Y] = \int_{x\in\Omega} dx\;\, x|\,\Psi(x,Y)\,|^2 \quad (8.7)$$

And it is trivial to calculate the *CP* of *X* given *Y*, starting from (6.3) or (8.2):

$$\boldsymbol{CP}\text{: } P(X\,|\,Y) = \frac{\mathrm{Tr}\left\{\rho_{X\cap Y}\right\}}{\mathrm{Tr}\left\{\rho_Y\right\}} = \frac{|\,\Psi(X,Y)\,|^2}{\int_{x\in\Omega} dx\,|\,\Psi(x,Y)\,|^2} \quad (8.8)$$

Our results (8.6-8.8) are expressed in terms of *probability density*, consistent with *Born statistical interpretation postulate* (see §2.8 of [6] or §1.2.1 of [7]). Moreover, they only use simple calculus, without reference to Measure theory or to W*-algebra.

Moreover, even if we are to deal with non-commutative observables, we may still use our unified expressions (6.1-6.3) to some extent, as we have demonstrated in §7.

Every standard textbook of introductory QM discusses expectation values of observables. But almost none of them even mention about conditional expectation. Why not? Actually, this is our original motivation to write this article. We wish that, with the help of our *simple* and *unified* expressions, the concept of *conditional expectation* now could be included in the context of introductory QM.

Hopefully, our unified expressions in (6.1-3) and (6.16-18) may also be useful for study of *conditional density operator* (*CDO*), *conditional probability* (*CP*) and *conditional expectation* (*CE*) in advanced quantum modeling [11-15].